\documentclass[11pt]{amsart}
\usepackage{cite}
\usepackage{mathrsfs} 
\usepackage{ifthen}
\usepackage{fancyhdr}
\usepackage{exscale}
\usepackage{tabularx}
\usepackage{latexsym}
\usepackage{amsmath,amssymb,amstext,amsthm} 
\usepackage{graphicx,color,epsfig}
\usepackage[table]{xcolor}
\usepackage{graphicx}
\usepackage{float}
\usepackage{verbatim}
\usepackage{enumerate}

\newtheorem{lem}{Lemma}
\newtheorem{thm}[lem]{Theorem}
\newtheorem{cor}[lem]{Corollary}

\def\hypergeom#1#2#3#4#5{{}_#1 F_{#2}\left({#3\atop#4}; \ #5\right)}

\begin{document}
\title{A Fourier analysis approach to disprove the weak Shanks conjecture}
\author{Jeffrey S. Geronimo}\address{School of Mathematics, Georgia Institute of Technology, 225 North Ave, Atlanta, GA 30332.} \email{jeffrey.geronimo@math.gatech.edu}

\author{Hugo J. Woerdeman}
\address{Department of Mathematics, Drexel University,
3141 Chestnut Street, Philadelphia, PA 19104}
\curraddr{}
\email{hugo@math.drexel.edu}
\thanks{HW is partially supported by National Science Foundation grant DMS 2348720.}


\maketitle

\begin{abstract}
We derive formulas for the Fourier coefficients of $|f|^2$, where $f(z_1,z_2)=(1-\frac{z_1+z_2}{r})^{-\alpha}$, in terms of hypergeometric functions. Using these formulas we provide additional counterexamples to the weak Shanks conjecture, which was recently disproven by B\'en\'eteau, Khavinson and Seco. The obtained formulas allow for (numerical) optimization over the parameters $\alpha$ and $r$. \end{abstract}



\section{Introduction}
\label{sec:intro}

Let $H^2 ({\mathbb D}^2)$ denote the Hardy Hilbert space on the bidisk ${\mathbb D}^2$, where ${\mathbb D} = \{ z \in {\mathbb C} : |z|<1 \}$. We let ${\mathcal P}_{n}$ denote all polynomials in two variables of total degree $\le n$. 
For each $f \in H^2 ({\mathbb D}^2)$ and
$ n \in {\mathbb N}$, a polynomial $p\in {\mathcal P}_{n}$ that minimizes the norm $\| 1-pf \|_{ H^2 ({\mathbb D}^2)}$ is called an {\em optimal polynomial approximant} (OPA)
for $\frac1f$ of degree $n$. 
After disproving the original Shanks conjecture following \cite{STJ}, Delsarte, Genin and Kamp [10] stated the 'weak Shanks conjecture': 

{\em Suppose $f$ is a polynomial with no zeros in ${\mathbb D}^2$, then its optimal polynomial 

approximants (OPAs) are zero-free in ${\mathbb D}^2$.}

\noindent In \cite{BKS} the authors recently provided a counterexample to the weak Shanks conjecture. In fact, they show that for the two variable function
$$ f(z_1,z_2)=(\frac{1}{1-\frac{z_1+z_2}{\sqrt{6}}})^{\frac52}$$ the optimal polynomial approximant $p_1(z_1,z_2)$ of degree 1 has zeroes in the bidisk ${\mathbb D}^2$, and consequently a counterexample to the weak Shanks conjecture is found by taking a take a sufficiently large degree Taylor polynomial of $f$. 
In this note we obtain formula's for the Fourier coefficients of $|f|^2$, where $f(z_1,z_2)=(1-\frac{z_1+z_2}{r})^{-\alpha}$, in terms of hypergeometric functions. This allows us to provide more counterexamples to the weak Shanks conjecture, as well as analyzing how far inside the bidisk the roots of the degree one OPA may be. 

\section{Fourier coefficients of squared modulus of a class of functions}

As we will see, it will be useful to have formulas for $|f|^2$, where $f \in H^2 ({\mathbb D}^2)$. For the functions appearing in \eqref{p2d}, we can express the Fourier coefficients in terms of hypergeometric functions. 
The hypergeometric functions we consider are defined for $|z|<1$  by the power series
$${}_{p}F_q\left(\begin{matrix}a_1,&\ldots&,a_{p}\\b_1,&\ldots&,b_q\end{matrix};\ z\right) = \sum_{n=0}^\infty \frac{(a_1)_n\cdots(a_{p})_n}{(b_1)_n\cdots(b_q)_n} \frac{z^n}{n!} .
$$
Here $(\alpha)_n$ is the (rising) Pochhammer symbol, which is defined by:
$$(\alpha)_n = \begin{cases}  1  & n = 0, \\
  \alpha (\alpha +1) \cdots (\alpha +n-1) & n > 0.
 \end{cases}$$
 \begin{thm}
Consider $f \in H^2 ({\mathbb D}^2)$ given by
\begin{equation}\label{p2d}
f(z_1,z_2)=(\frac{1}{1-\frac{z_1+z_2}{r}})^{\alpha},
\end{equation}
where $r>2$ and $\alpha \ge 0$. Denote the Fourier coefficients of $|f|^2$ by $ c_{k_1,k_2}$, $k_1,k_2 \in {\mathbb Z}$. 
Then, for $k_1,k_2\ge 0$, we have 
\begin{align}\nonumber
c_{k_1,k_2}&=\frac{1}{r^{k_1+k_2}}\frac{(\alpha)_{k_1+k_2}}{k_1!k_2!}\hypergeom{4}{3}{\alpha,\alpha+k_1+k_2,\frac{k_1+k_2}{2}+1,\frac{k_1+k_2+1}{2}}{k_1+k_2+1,k_1+1,k_2+1}{\frac{4}{r^2}}
\end{align}
and
\begin{align}\nonumber
c_{k_1,-k_2}=\frac{1}{r^{k_1+k_2}}\frac{(\alpha)_{k_1}(\alpha)_{k_2}}{k_1!k_2!}\hypergeom{4}{3}{\alpha+k_1,\alpha+k_2,\frac{k_1+k_2+1}{2},\frac{k_1+k_2}{2}+1}{k_1+1,k_2+1,k_1+k_2+1}{\frac{4}{r^2}}.
\end{align}
\end{thm}

We note that for the case when $\alpha=1$ the above formulas for the Fourier coefficients were derived in \cite{GWW}.

\begin{proof} First observe that
\begin{align}\label{2dinv}
f(z_1,z_2)&=\sum_{i=0}^{\infty}\frac{(\alpha)_i}{i!r^i}(z_1+z_2)^i\nonumber\\&=\sum_{i=0}^{\infty}\frac{(\alpha)_i}{i!r^i}\sum_{j=0}^i\frac{i!}{j!(i-j)!}z_1^{j}z_2^{i-j}.
\end{align}
Therefore
\begin{align}\label{3dsp}
& |f(z_1,z_2)|^2=\nonumber\\ 
&=\sum_{i_1=0}^{\infty}\frac{1}{r^{i_1}}\sum_{j_1=0}^{i_1}\frac{(\alpha)_{i_1}}{j_1!(i_1-j_1)!}z_1^{j_1}z_2^{i_1-j_1}\sum_{i_2=0}^{\infty}\frac{1}{r^{i_2}}\sum_{j_2=0}^{i_2}\frac{(\alpha)_{i_2}}{j_2!(i_2-j_2)!}z_1^{-j_2}z_2^{-i_2+j_2}\nonumber\\&=:\sum_{k_1=-\infty}^{\infty}\sum_{k_2=-\infty}^{\infty}c_{k_1,k_2}z_1^{k_1}z_2^{k_2}.
\end{align}
We start with $k_1\ge0$ and $k_2\ge0$ and  matching coefficients we find $k_2=j_1-j_2$ and $k_1=i_1-i_2-j_1+j_2$. Rewriting as $j_1=k_2+j_2$ and $i_1=k_1+k_2+i_2$ we find (with $i=i_2$ and $j=j_2$),
$$
c_{k_1,k_2}=\sum_{i=0}^{\infty}\frac{(\alpha)_i(\alpha)_{k_1+k_2+i}}{r^{2i+k_1+k_2}}\sum_{j=0}^{i}\frac{1}{j!(i-j)!}\frac{1}{(k_1+i-j)!(k_2+j)!}.
$$
Since 
$$
(k_2+j)!=(k_2+1)_j k_2!,
$$
$$
\frac{1}{(i-j)!}=\frac{(-1)^j(-i)_j}{i!},
$$
and
$$
(k_1+i-j)!=(-1)^j\frac{(k_1+i)!}{(-k_1-i)_j},
$$
we find, using the Chu-Vandermonde formula \cite[p.~67]{aar}, that the sum on $j$ can be written as
$$
\sum_{j=0}^i\frac{(-i)_j(-k_1-i)_j}{(1)_j(k_2+1)_j}=\frac{(k_1+k_2+i+1)_i}{(k_2+1)_i}.
$$
With the identities 
$$
\frac{(k_1+k_2+i)!}{(k_1+i)!}=\frac{(k_1+k_2)!(k_1+k_2+1)_i}{k_1!(k_1+1)_i},
$$
$$
(\alpha)_{k_1+k_2+i}=(\alpha)_{k_1+k_2}(\alpha+k_1+k_2)_i ,
$$
and
$$
(k_1+k_2+i+1)_i=4^i\frac{(\frac{k_1+k_2}{2}+1)_i(\frac{k_1+k_2+1}{2})_i}{(k_1+k_2+1)_i },
$$
we have
\begin{align}\label{ck1k2}
c_{k_1,k_2}&=\frac{1}{r^{k_1+k_2}}\frac{(\alpha)_{k_1+k_2}}{k_1!k_2!}\sum_{i=0}^{\infty}\frac{4^i}{r^{2i}}\frac{(\alpha)_i(\alpha+k_1+k_2)_i(\frac{k_1+k_2}{2}+1)_i(\frac{k_1+k_2+1}{2})_i}{(1)_i(k_1+k_2+1)_i(k_2+1)_i(k_1+1)_i}\nonumber\\&=\frac{1}{r^{k_1+k_2}}\frac{(\alpha)_{k_1+k_2}}{k_1!k_2!}\hypergeom{4}{3}{\alpha,\alpha+k_1+k_2,\frac{k_1+k_2}{2}+1,\frac{k_1+k_2+1}{2}}{k_1+k_2+1,k_1+1,k_2+1}{\frac{4}{r^2}}.
\end{align}

Next  we compute $c_{k_1,-k_2}$ with $k_1, k_2\ge0$. Using equation~\eqref{3dsp}  restricted to this region and setting $i_1=k_1+i_2$ and $j_1=j_2-k_2$  gives
\begin{align*}
&\sum_{k_1=0}^{\infty}\sum_{k_2=0}^{\infty}c_{k_1,-k_2}z_1^{k_1}z_2^{-k_2}\\&=\sum_{k_1=0}^{\infty}\sum_{j_2=0}^{\infty}\sum_{k_2=0}^{j_2}\sum_{i_2=0}^{\infty}\frac{z_1^{k_1}z_2^{-k_2}}{r^{k_2-k_1+2i_2+2j_2}}\frac{(\alpha)_{k_1+i_2+j_2-k_2}(\alpha)_{i_2+j_2}}{(j_2-k_2)!(i_2+k_1)!j_2!i_2!}
\end{align*}
Equating coefficients in the above equation then setting $i=i_2+j_2-k_2$ and $j=j_2-k_2$ yields
\begin{equation*}
c_{k_1,-k_2}=\sum_{i=0}^{\infty}\frac{1}{r^{2i+k_1+k_2}}\sum_{j=0}^i\frac{(\alpha)_{i+k_1}(\alpha)_{i+k_2}}{j!(i+k_1-j)!(j+k_2)!(i-j)!}.
\end{equation*}
Substitution of the above formulas into the sum on $j$ gives
\begin{align*}
\sum_{j=0}^i\frac{1}{j!(i+k_1-j)!(j+k_2)!(i-j)!}&=\frac{1}{i!k_2!(i+k_1)!}\sum_{j=0}^i\frac{(-i)_j(-k_1-i)_j}{j!(k_2+1)_j}\\&=\frac{1}{i!k_2!(i+k_1)!}\frac{(k_1+k_2+i+1)_i}{(k_2+1)_i}\\&=\frac{4^i}{i!(i+k_1)!}\frac{(\frac{k_1+k_2}{2}+1)_i(\frac{k_1+k_2+1}{2})_i}{(k_1+k_2+1)_i},
\end{align*}
where the Chu-Vandermonde formula  was used to obtain the second equality above.
Thus after some simplification we find
\begin{align}\label{ck2negf}
c_{k_1,-k_2}=\frac{1}{r^{k_1+k_2}}\frac{(\alpha)_{k_1}(\alpha)_{k_2}}{k_1!k_2!}\hypergeom{4}{3}{\alpha+k_1,\alpha+k_2,\frac{k_1+k_2+1}{2},\frac{k_1+k_2}{2}+1}{k_1+1,k_2+1,k_1+k_2+1}{\frac{4}{r^2}}.
\end{align}
\end{proof}

\section{Optimal polynomial approximants}

The equation \eqref{YW} below provides a way to compute the optimal polynomial approximants. This result has appeared earlier; see, e.g., \cite[Theorem 2.1]{BC}. We use the multivariable notation $z^k = z_1^{k_1}z_2^{k_2}$, where $z=(z_1,z_2)$ and $k=(k_1,k_2)$. In addition ${\mathbb N}=\{ 1,2, \ldots \}$ and ${\mathbb N}_0 = {\mathbb N}\cup \{ 0 \}$.

\begin{thm} Given $f=\sum_{k\in {\mathbb N}_0^2} f_kz^k \in H^2 ({\mathbb D}^2)$. Denote the Fourier coefficients of $|f|^2$ by $c_k, k\in{\mathbb Z}$. Let $n\in {\mathbb N}$, $\Lambda_n = \{ (k_1,k_1) \in{\mathbb N}_0^2 : k_1+k_2\le n \}$ and $p(z) = \sum_{k \in \Lambda_n} p_kz^k$ be the optimal polynomial approximant of $\frac1f$ of degree $n$. Then
\begin{equation}\label{YW} (c_{k-l})_{k,l\in \Lambda_n} (p_l)_{l\in\Lambda_n} = (\overline{f_k}\delta_{0,k})_{k\in\Lambda_n}, \end{equation}
where $\delta_{0,k}=1$ when $k=0$ and $\delta_{0,k}=0$ otherwise.
\end{thm}  

To write out equation \eqref{YW}, one needs to order the elements of $\Lambda_n$. If, for instance, we take $n=1$ and order $\Lambda_1 = \{ (0,0), (0,1), (1,0)\}$ lexicographically, \eqref{YW} becomes 
\begin{equation}\label{1} \begin{pmatrix} c_{00} & c_{0,-1} & c_{-1,0} \cr c_{0,1} & c_{00} & c_{-1,1} \cr c_{1,0} & c_{1,-1} & c_{00} \end{pmatrix}  \begin{pmatrix} p_{00}  \cr p_{01}  \cr p_{10} \end{pmatrix}  = \begin{pmatrix} \overline{f_{00}}  \cr 0  \cr 0 \end{pmatrix} . \end{equation}

\begin{proof} Notice that minimizing the norm $\| 1-pf \|_{ H^2 ({\mathbb D}^2)}$ comes down to computing the distance from 1 to the linear space ${\mathcal P}_n f$. Thus at the optimum $p$ we must have that $1-pf \perp {\mathcal P}_n f$. Since ${\mathcal P}_n f = {\rm span} \{ z^k f : k \in \Lambda_n \}$, the orthogonality relation comes down to the equations
\begin{equation}\label{2} \langle 1-pf , z^kf \rangle_{H^2 ({\mathbb D}^2)} = 0, k \in \Lambda_n. \end{equation}
Next, we observe that 
$$ \langle 1 , f \rangle_{H^2 ({\mathbb D}^2)} = \overline{f_0}, \ \langle 1 , z^kf \rangle_{H^2 ({\mathbb D}^2)} = 0, k\in \Lambda_n\setminus \{ 0 \},$$ and
$$  \langle pf , z^kf \rangle_{H^2 ({\mathbb D}^2)} = \sum_{l\in\Lambda_n} p_l \langle z^{l-k} , |f|^2 \rangle_{L^2 ((\partial{\mathbb D})^2)} = \sum_{l\in\Lambda_n} p_l \overline{c_{l-k}} = \sum_{l\in\Lambda_n} p_l c_{k-l}. $$
Thus we see that \eqref{2} corresponds exactly to the matrix equation \eqref{YW}.
\end{proof}

If we apply the above to the function \eqref{p2d} and let $n=1$, we obtain
$$  a:=c_{00}=\hypergeom{4}{3}{\alpha,\alpha,\frac{1}{2},1}{1,1,1}{\frac{4}{r^2}}=
\hypergeom{3}{2}{\alpha,\alpha,\frac{1}{2}}{1,1}{\frac{4}{r^2}}, $$
$$ b:=c_{0,\pm1}=c_{\pm1,0} =\frac{\alpha}{r}\hypergeom{4}{3}{\alpha,\alpha+1,\frac{3}{2},1}{2,1,2}{\frac{4}{r^2}} =\frac{\alpha}{r}
\ \hypergeom{3}{2}{\alpha,\alpha+1,\frac{3}{2}}{2,2}{\frac{4}{r^2}},
$$ and
$$ c:= c_{-1,1}=c_{1,-1}= \frac{\alpha^2}{r^2}\ \hypergeom{4}{3}{\alpha+1,\alpha+1,2,\frac{3}{2}}{2,2,3}{\frac{4}{r^2}} $$ $$\ \ \ \ \ \ \ \ \ \ \ \ \ \ \ \ \ \ \ \ \ \ =\frac{\alpha^2}{r^2}\ \hypergeom{3}{2}{\alpha+1,\alpha+1,\frac{3}{2}}{2,3}{\frac{4}{r^2}}.$$
Note that for $f$ as in \eqref{p2d}, we have that $f(0,0)=1$. In addition, $f$ is symmetric in the variables $z_1,z_2$, so the same holds for its OPAs. Thus the OPA for $\frac1f$ of degree $1$ equals $p(z_1,z_2)=p_0 + p_1(z_1+z_2)$,
where 
$$ \begin{pmatrix} p_0 \cr p_1 \cr p_1 \end{pmatrix} = 
\begin{pmatrix} a & b & b \cr b & a & c \cr b & c & a \end{pmatrix}^{-1} \begin{pmatrix} 1 \cr 0 \cr 0 \end{pmatrix}. $$ Note that $p$ has a root in the open bidisk exactly when $2|p_1|>|p_0|$. Computing its fraction, 
we find that 
$$ \frac{|p_0|}{2|p_1|}= \frac{|a^2-c^2|}{2 |b(a-c)|} = \frac{a+c}{2b}, $$ where we used that $a,b,c>0$.
Thus we have the following corollary.

\begin{cor}\label{cor} Let $f$ be as in \eqref{p2d}. The following are equivalent.
\begin{itemize}
\item[(i)] The OPA for $\frac1f$ of degree $1$ has a root in the open bidisk.
\item[(ii)]
\begin{equation*}\label{ineq} 2 \frac{\alpha}{r} \ \hypergeom{3}{2}{\alpha,\alpha+1,\frac{3}{2}}{2,2}{\frac{4}{r^2}} > \hypergeom{3}{2}{\alpha,\alpha,\frac{1}{2}}{1,1}{\frac{4}{r^2}} +\frac{\alpha^2}{r^2}\ \hypergeom{3}{2}{\alpha+1,\alpha+1,\frac{3}{2}}{2,3}{\frac{4}{r^2}}.\end{equation*}
\item[(iii)] \begin{align*}
&\int_0^{\infty}\int_0^{\infty}\int_0^1e^{-t-y}t^{\alpha-1}y^{\alpha-1}(\frac{4}{r}yw^{1/2}+w^{-1/2})(1-w)^{-1/2} \hypergeom{0}{1}{-}{2}{\frac{4}{r^2}tyw}{\rm d}t{\rm d}y{\rm d}w\\& >2\int_0^{\infty}\int_0^{\infty}\int_0^1 e^{-t-y}t^{\alpha-1}y^{\alpha-1} w^{-1/2}(1-w)^{-1/2} \hypergeom{0}{1}{-}{1}{\frac{4}{r^2}tyw} {\rm d}t{\rm d}y{\rm d}w.
\end{align*}
\end{itemize}
\end{cor}

The potential advantage of item (iii) above may be that the dependence on $\alpha$ is more conducive to further analysis.

\begin{proof} We proved the equivalence of (i) and (ii) before the statement of the corollary. It remains to prove the equivalence of (ii) and (iii). 
Observe that \cite[p. 115]{aar}
%
$$
\hypergeom{3}{2}{a_1,a_2,a_3}{b_1,b_2}{x}=\frac{1}{\Gamma(a_1)\Gamma(a_2)}\int_0^{\infty}\int_0^{\infty} e^{-t-y} t^{a_1-1}y^{a_2-1}\hypergeom{1}{2}{a_3}{b_1,b_2}{xyt}dtdy,
$$
where $a_3\in\{\frac12 , \frac32 \}$ and the $b_1,b_2\in\{1,2,3\}$. Plugging this in (ii) leads to 
\begin{align*}
&2 \frac{\alpha}{r} \frac{1}{\Gamma(\alpha)\Gamma(\alpha+1)}\int_0^{\infty}\int_0^{\infty}e^{-t-y}t^{\alpha-1}y^{\alpha} \hypergeom{1}{2}{\frac{3}{2}}{2,2}{\frac{4}{r^2}ty}\\& >\frac{1}{\Gamma(\alpha)\Gamma(\alpha)}\int_0^{\infty}\int_0^{\infty}e^{-t-y}t^{\alpha-1}y^{\alpha-1} \hypergeom{1}{2}{\frac{1}{2}}{1,1}{\frac{4}{r^2}ty}\\& +\frac{\alpha^2}{r^2}\frac{1}{\Gamma(\alpha+1)\Gamma(\alpha+1)}\int_0^{\infty}\int_0^{\infty}e^{-t-y}t^{\alpha}y^{\alpha} \hypergeom{1}{2}{\frac{3}{2}}{2,3}{\frac{4}{r^2}ty}.
\end{align*}
Since $\Gamma(\alpha+1)=\alpha\Gamma(\alpha)$ and $\alpha>0$ the above inequality becomes
\
\begin{align*}
&\frac{2}{r} \int_0^{\infty}\int_0^{\infty}dtdye^{-t-y}t^{\alpha-1}y^{\alpha} \hypergeom{1}{2}{\frac{3}{2}}{2,2}{\frac{4}{r^2}ty}\\& >\int_0^{\infty}\int_0^{\infty}dtdye^{-t-y}t^{\alpha-1}y^{\alpha-1} \hypergeom{1}{2}{\frac{1}{2}}{1,1}{\frac{4}{r^2}ty}\\& +\frac{1}{r^2}\int_0^{\infty}\int_0^{\infty}dtdye^{-t-y}t^{\alpha}y^{\alpha} \hypergeom{1}{2}{\frac{3}{2}}{2,3}{\frac{4}{r^2}ty}.
\end{align*}
Substitute the relation
$$
\hypergeom{1}{2}{a}{b,c}{z}=\frac{\Gamma(b)}{\Gamma(a)\Gamma(b-a)}\int_0^1 t^{a-1}(1-t)^{b-a-1}\hypergeom{0}{1}{-}{c}{zt}dt
$$
into the above inequality to find (using $\Gamma(\frac32 )=\frac12 \Gamma(\frac12 )$)
\begin{align*}
&\frac{4}{r} \int_0^{\infty}\int_0^{\infty}\int_0^1dtdydwe^{-t-y}t^{\alpha-1}y^{\alpha}w^{1/2}(1-w)^{-1/2} \hypergeom{0}{1}{-}{2}{\frac{4}{r^2}tyw}\\& >\int_0^{\infty}\int_0^{\infty}\int_0^1dtdydwe^{-t-y}t^{\alpha-1}y^{\alpha-1} w^{-1/2}(1-w)^{-1/2} \hypergeom{0}{1}{-}{1}{\frac{4}{r^2}tyw}\\& +\frac{2}{r^2}\int_0^{\infty}\int_0^{\infty}\int_0^1dtdydwe^{-t-y}t^{\alpha}y^{\alpha} w^{1/2}(1-w)^{-1/2} \hypergeom{0}{1}{-}{3}{\frac{4}{r^2}tyw}.
\end{align*}
Finally, use the recurrence
$$
\hypergeom{0}{1}{-}{1}{z}-\hypergeom{0}{1}{-}{2}{z}=\frac{z}{2}\ \hypergeom{0}{1}{-}{3}{z}
$$ to arrive at the inequality in (iii).
\end{proof}

Applying Corollary \ref{cor}(ii) for the choice to $\alpha =\frac52$ and $ r=\sqrt{6}$, gives
$$ 84.4198= 2 \cdot 42.2099 > 45.3768 + 37.0258 = 82.4023, $$
confirming the counterexample by \cite{BKS}.

The above test provides a way to produce many counterexamples, including, for example, by taking high enough degree Taylor polynomials of the rational functions
$$ \left( 1- \frac{z_1+z_2}{2.1} \right)^{-2} , \left( 1- \frac{z_1+z_2}{3} \right)^{-3}, $$
for which Corollary \ref{cor}(ii) holds as
$$ 558.8698= 2\cdot 279.4349 > 285.4070  + 266.1290 = 551.5360 $$
and
$$ 28.7862= 2 \cdot 14.3931 > 16.5183 +  11.6884 = 28.2067, $$
respectively. In fact, the Matlab graph below shows the region of pairs of $(\alpha , r )$ for which Corollary 3(ii) holds.

\vspace{-1.252cm}\hspace{0.2cm}\includegraphics[height=8cm, width=12cm]{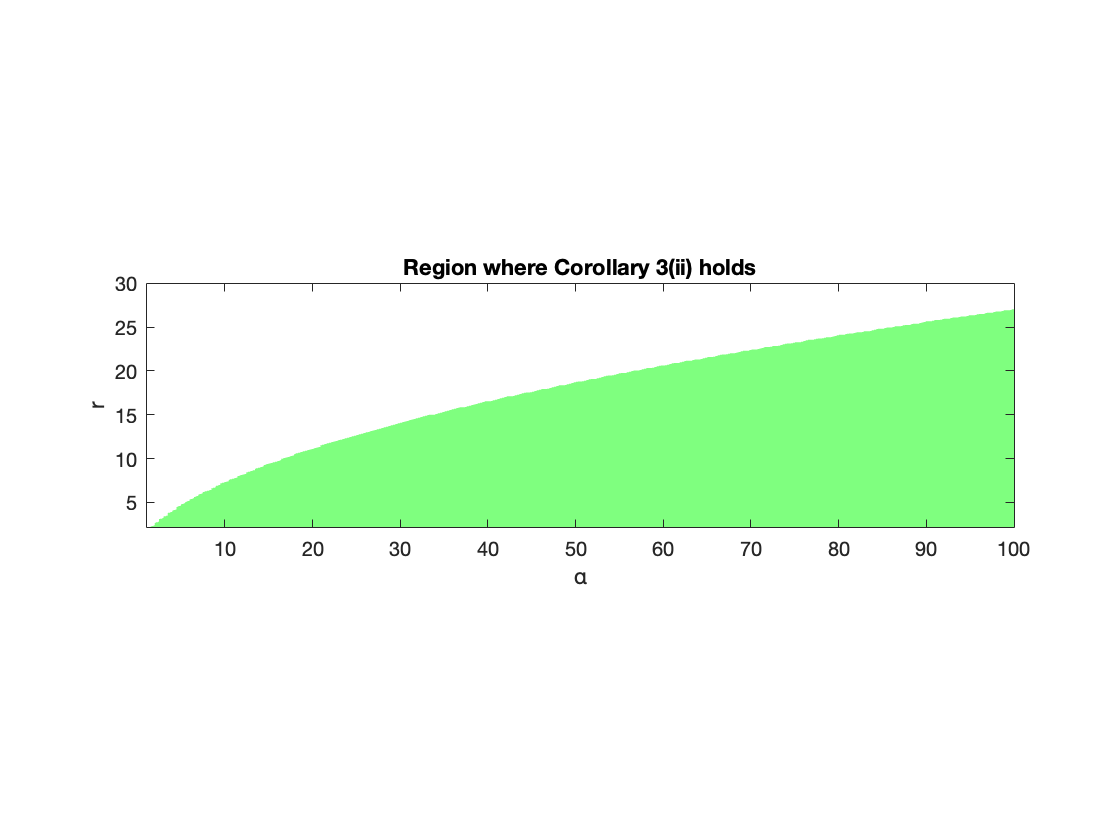}

\vspace{-2.4cm}

In the last section of \cite{BKS} the question was asked how close to 0 the roots of a degree one OPA of $\frac1f$ can be. To address this question we used the Matlab command fmincon to minimize the quantity \begin{equation}\label{quotient} \frac{\hypergeom{3}{2}{\alpha,\alpha,\frac{1}{2}}{1,1}{\frac{4}{r^2}} +\frac{\alpha^2}{r^2}\ \hypergeom{3}{2}{\alpha+1,\alpha+1,\frac{3}{2}}{2,3}{\frac{4}{r^2}}}{2 \frac{\alpha}{r} \ \hypergeom{3}{2}{\alpha,\alpha+1,\frac{3}{2}}{2,2}{\frac{4}{r^2}}}.\end{equation} We find the minimum at the values
$r=2.533672086469380$, $\alpha=
2.551918826591946 $ with the quotient equalling 
$0.975766335259681$; see also the Mathematica graph below. For this value of $r$ and $\alpha$ we find the degree one OPA of $\frac1f$ to equal
$$0.5018   -0.2571(z_1+z_2), $$
which, for instance, has a root at $z_1=z_2= 0.975766335259681$. This OPA has roots that are the closest (in $\| \cdot \|_1$, for instance) to the origin.

\hspace{0.2cm}\includegraphics[height=7.5cm, width=10cm]{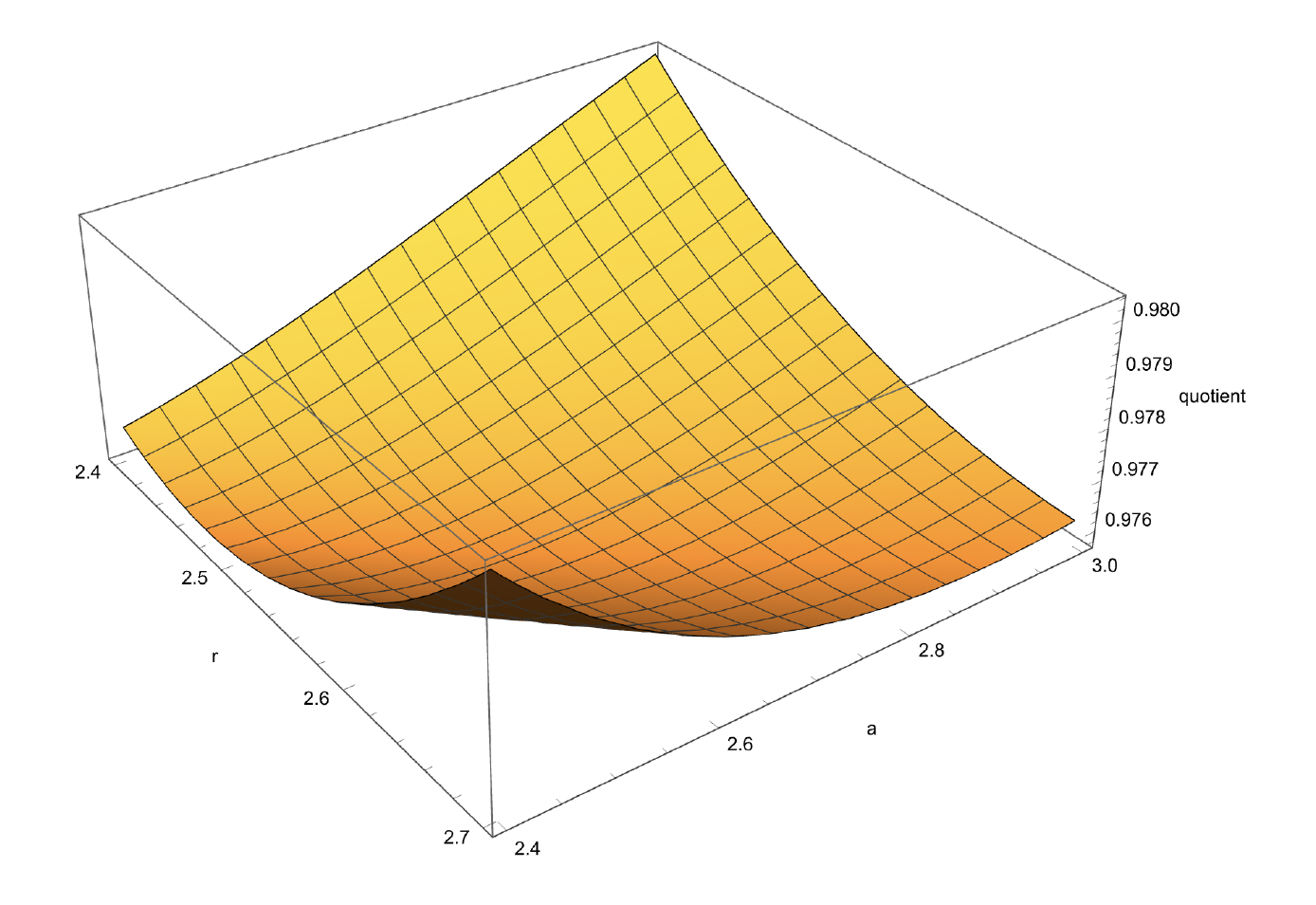}

\centerline{Graph of the quotient \eqref{quotient} where $r\in [2.4, 2.7]$ and $\alpha \in[ 2.4, 3]$.}

\bibliographystyle{plain}

\end{document}